\documentclass[twoside,letterpaper,draft,12pt]{amsart}

\usepackage{amsmath,amsthm,latexsym,xspace,amscd,amssymb, textcomp}
\usepackage[margin=1.3 in]{geometry}

\title{ALGEBRAS WITH LAURENT POLYNOMIAL IDENTITY }

\title{ALGEBRAS WITH LAURENT POLYNOMIAL IDENTITY}    
\author{Claudenir Freire, Ramon Codamo}
\address{Departmento de Matemática, Universidade Federal do Amazonas, Amazonas, Manaus, Brasil}
\email{claudiomao1@gmail.com}
\email{ramoncodamo@gmail.com}

\newtheorem{teo1}{Theorem}[section]

\newtheorem{lem1}[teo1]{Lemma}
\newtheorem{cor1}[teo1]{Corollary}
\newtheorem{prop}[teo1]{Proposition}

\theoremstyle{remark}
\theoremstyle{remark}

\def\ndv{\ {\mid \kern -0.7 em {\scriptstyle \not}} \ \ }

\def\nd{\ {\mid \kern -0.4 em {\scriptstyle \not}} \ \ }

\begin{document}


\thispagestyle{empty}
\begin{titlepage}
	\begin{center} 
		\vspace{0.7cm}
		\begin{center}
			
			
			
			\vspace{4cm}
			
			\begin{center}
				{ \bf ALGEBRAS WITH LAURENT POLYNOMIAL IDENTITY}\\
				\vspace{1cm}                    
				By
			\end{center}

			{\normalsize \vspace{1cm} }
			
			
			{\normalsize \setcounter{footnote}{1} \vspace{10mm} \renewcommand{%
					\thefootnote}{\fnsymbol{footnote}} {\large \textbf{Claudenir
						Freire Rodrigues, Ramon C\'odamo B. da Costa }}}
			

			
			\vspace{2cm} \setcounter{footnote}{1}
			\renewcommand{\thefootnote}{ \fnsymbol{footnote}}
			
			
			\vspace{2cm}
			
		\end{center}
		\vspace{0.5cm}
	\end{center}
	\vspace{8cm}
	{\footnotesize Universidade Federal do Amazonas - Departamento de Matem\'atica\\
		Av. General Rodrigo Octavio Jord\~ao Ramos, 1200 - Coroado I, Manaus - AM - Brasil, 69067-005}\\
	{\footnotesize claudiomao1@gmail.com (corresponding author)}\\
	{\footnotesize ramoncodamo@gmail.com}
	
\end{titlepage}

\newpage

\begin{abstract}
In this article we shows some results about algebra with the group of units having special polynomial identity.
\end{abstract}
\maketitle

\section{\textbf{Introduction}}

A Laurent polynomial $ f = f(x_1,...,x_l) $ in the noncommutative variables $x_i, \, i=1,\ldots,l$ is an element non-zero in the group algebra  $RF_l$ over a ring  $R$  with free group $F_l = < x_1,...,x_l >$. One says  $ f $ in $RF_l $ is a Laurent Polynomial Identity  ($LPI$) for an $R$ algebra $A$ (respc. for $U(A)$ the group of units in A) if $f(a_1, ..., a_l)=0 $ for all sequence $a_1,..., a_l$ in $A$ (respc. in U(A)). A word  $ w $ in $RF_l$ is group identity for all $U(A)$ if $w(a_1, ..., a_l)=1$, for  $a_1,..., a_l$ in $U(A)$ .\\ 

\hspace*{0,05cm}\textbf{Brian Hartley´s Conjecture}:	Let G be a torsion group and R a field. If the unit
group $U(RG)$ of $RG$ satisfies a group identity $w=1$, then RG satisfies a polynomial
identity.\\

In [2] A.Giambruno, E.Jespers and A.Valenti shows that the Brian Hartley´s Conjecture is true for a group algebra $RG$ over an infinite commutative domnain $R$ and a torsion group $G$ and $G$ has no divisible order elements by $p$ with characteristic of $R$ is $p$. For this purpose they proved the following crucial result.

\begin{prop} 
	
	Let $A$ be an algebra  over an infinite commutative domain $R$ and suppose that $U(A)$ satisfies a group identity. There exists a positive m such that if $ \, a,b,c, u $ in $A$ and $a^{2}=bc=0$, then $bacA$ is nil right ideal of bounded exponent less or equal than m.
	
\end{prop}	
	  
The case identity group  $w=1$ in the conjecture says  the  unity group satisfy equivalently the Laurent polynomial  $F=  1- w  $.  We use this idea and study the case  $F=\, a_1 + a_2w_2 +...+a_n w_n $. Clearly the case $F=  1- w  $ is a special  case of our generalization.

In this paper we prove the following results. \\
 
 \newpage
 
\begin{prop}\label{p1}
	Let $A$ be an algebra with the group of units $U(A)$ admits a  $LPI$  over a ring R whith unit  whose non-constants words has the sum of exponents  non-zero at least one of the variables. Then there exists a polynomial $f \, \in \, R[X]$ with the limited degree $d\leq 4(-l+r)+3$, where $l=min \, \{\sum expw_i\}$ and $r=max\{\sum exp \, w_i \}$, determinated  by the $LPI$  such that for all $ \, a,b,c, u $ in $A$ with $a^{2}=bc=0$,  $f(bacu)=0$. In particular, $bacA$  is an algebraic ideal.
	
\end{prop}
As a consequence of Proposition 1.2 we prove the following: 

\begin{cor1}\label{c1} 
	Let $A$ be an algebra over a commutative domain R with $ \mid R \mid > d  $ whose unit group $U(A)$ has a $LPI$. If $a,\, b,\, c \in A$ with $a^2=bc=0 $, then $bacA$ is a nil ideal with limited exponent.	
	
\end{cor1}

Notice that the restriction on $LPI$ is necessary because we will show that without this assumption Proposition 1.2 i not valid. In this article we always adopt $LPI$ with this restriction. Also, by changin $x_i$ by $x^{-i}yx^i$ we may assume that  $LPI$ in two variables.

\newpage


	
%
	

\textbf{\textit{Proof of Proposition \ref{p1}}: }

 Let  $P\,=\, a_1 + a_2w_2 +...+ a_r w_t$ the $LPI$ with $w_i $ in the form     $$w_i=x_1^{r_1}x_2^{s_2} \cdots x_1^{r_k}x_2^{s_k}, k\geq 1, r_i,\, s_i $$ integers. If $ \sum exp \, w_i = \sum exp _{x_1}\, w_i + \sum exp _{x_2}\, w_i = 0 $ then we  substitute $ x_1= x_1^k$ or $x_2= x_2^k$ with $k > 1$ big enough,
	We get a new $LPI$ in the form  $P\,=\, a_1 + a_2{w'}_2 +...+ a_r {w'}_r,$ where   
	$$ \sum exp \, {w'}_i = \sum exp _{x_1}\, {w'}_i + \sum exp _{x_2}\, {w'}_i = k \sum exp _{x_1}\, w_i + \sum exp _{x_2}\, w_i $$  or $$\sum exp \, {w'}_i = \sum exp _{x_1}\, w_i +  k\sum exp _{x_2}\, w_i. $$     In any  case that is not zero.

	In  general,  $P(1,1)\,=\, a_1 + a_2 +...+ a_r =0$, so $ a_2 +...+ a_r \neq 0 $ , because $a_1 \neq 0$. Hence any collection of sums with distinct parcels involving all coefficients  $ a_2,..., a_r $  has a non-zero term.
	
	From this $$P(\alpha, \alpha)\,=\, a_1 + a_2 \alpha ^{\sum exp w_2} +...+ a_k \alpha ^{\sum exp w_k} =  0, \,  \, 1<k\leq r,  $$ is a non-zero polynomial with integers exponents  over $R$, for all $\alpha \in U(A)$. In the following we see that all exponents are not  necessarily positive.
	
	Let $\alpha ^{-l}$ be where $l=min \, \{\sum expw_i\}$, so $\alpha ^{-l}P(\alpha,\alpha)= f_0(\alpha)=0$, where $f_0\neq 0$ is a polynomial	with degree $ \leq -l + r$, over $R$. Therefore  $U(A)$ is algebraic over $R$. Now, we are going to use this polynomial with a particular unit in $U(A)$.
	First let $b=c$ be, so $a^2=b^2 =0$.    
	Then we have $(1 + aua), \, (1 + bauab) \in U(A)$, for all $u \in A$.
	Let $Q_n$ be the set of all  products in $aua$ and $bauab$ with at most $n-1$ factors, and let 
	$W$  be the set of all  products in coeficients $ a_i`s$, $ aua$  and $ bauab$.

	Thus $\alpha = (1 + aua).(1 + bauab) \, \in \, U(A) $, and for this unit,    
	$$\alpha ^n= 1 + aua + bauab +  \sum_{t\in Q_n}t+  (auabauab)^n$$ for all $n \geq 1$.

	It follows that there exists  a polynomial $f_0\neq 0 $ over $R$ such that
	
\begin{eqnarray*}
f_0(\alpha) &=&  a_1 + a_2 +...+ a_r + \sum_{t\in W}t + \hspace*{0,8cm}  \sum ^r_{i=2}  a_i (auabauab)^{s_i}\\&=&\sum_{t\in W}t  +\sum ^r_{i=2} \,  a_i (auabauab) ^{s_i} =0. 
\end{eqnarray*}

	As $ a_2 +...+ a_r \neq 0 $, after to organize  the   powers of $auabauab$  which appears in the last sum,  has a non-zero coefficients.  This term has the form $(a_{j_1}+...+a_{j_k})(auabauab)^s$ where $s \, \leq \, (-l + r)$. So, $abf_0(\alpha)au=0$ is a  non-zero polynomial $g$ over R such that $g(abau) = 0$ with degree $ \leq \,2( -l + r ) +1$. This finishs the  case  $b=c$.

	In general, if $a,\, \alpha ,\, \beta \in A$ and $a^{2}=\alpha \beta =0$, then  we have $(\beta u \alpha)^{2}=0 \,,$ 	for all $u \in A $. So $\beta u \alpha a \beta u \alpha A$ is algebraic and there exists $g$ over  $R$ which is not zero and
	
	$g( (\beta u \alpha a \beta u \alpha)a)=0$, from this there exists a non-zero polynomial $f=(\alpha a ) f_0 ( \beta u )$ with $f(\alpha a \beta u)=0$ such that $d=degree\,f \leq \, \, 4( -l + r) +3  $. Therefore $\alpha a \beta A$ is algebraic over $R$.

	
	

\newpage

\textbf{\textit {Proof of Corollary \ref{c1}} :}

	By Proposition \ref{p1}, there exists $f \in R[x]\setminus\{0\}$ with  degree  $ d \leq m = 4( -l + r ) +3$, such that $f(bac \lambda u)=0$ for all $\lambda \in R,\, u\, \in \, A$. That is, $\lambda p_1 + \lambda ^2 p_2 +...+ \lambda ^d  p_d=0$ where $p_i$ is polynomial on $ bacu $ with $p_d=b_d (bacu)^d$ . Now, by the Vandermond argument, $p_1=...=p_d=0$, so $ (bacu)^d=0$ .

The next results are LPI versions of Lema 3.2,  Lema 3.3  of \cite{LIU} and  Corolary 3  of  \cite{GJV94}. Their proofs are analogous to the  group identity case changing GI by LPI.

\begin{cor1}
	
	Let $A$ be a semiprime algebra over an  commutative domain $R$ with $|R|>d$. If $U(A)$ satisfies a $LPI$, then every idempotent element of $R^{-1} A$ is central.			
	
\end{cor1}	

\begin{cor1}\label{c2}
	
	Let $A$ be an algebra over a field $K$ with $ U(A)$ satisfying a $LPI$, and   $a,\, b \, \in A $ such that $a^2=b^2=0$.
	
	\begin{enumerate}
		
		\item If $\mid K \mid > 2d$, then $(ab)^{2d}=0.$ 
		\item If $ab$ is nilpotent, then $(ab)^{2d}=0.$
	\end{enumerate}

\end{cor1}

\begin{cor1}
	
	Let $K$ be any field and $n\geq 2$. Then the following are equivalent
	
	\begin{enumerate}
		\item  $U(M_n(K))$ satisfies the $LPI$. 
		\item  $K$ is a finite field.
		\item  $U(M_n(K))$ satisfies a group identity.
	\end{enumerate}
	
In addition, if one of these conditions is hold, then

	$\mid K \mid  \leq 2d $  and  
	$ n \leq 2log_{\mid K \mid}2d + 2 \leq 2log_2 2d + 2 . $

\end{cor1}

In  general,  Proposition 1.2 is not hold without the restrictions  about exponents. If it does, then $ bacu$ is algebraic for all $u \in A$, in particular, for $ c=b $ and $u=a$, $ba$ is algebraic.  It follows that by    Vandermonde augument with $R$ a infinite commutative domain $ba$ is nilpotent. In general,   $ba$ is not a nilpotent element. For example, let $A=M_n(R), \, n>1$. By $Amitsur-Levitzki$ theorem \cite{P} $A$ satisfies a $LPI$ 

$$ f_1(x_1,...,x_{2n})= S_{2n}\cdot(x_1\cdots x_{2n})^{-1}  $$  where $s_{2n}$ is the polynomial standard.  
However $A$ with $a=e_{21}$ and $b=e_{12}$,  $ba=e_{11}$ is not nilpotent. Note that this algebra also satisfies
$$f_2(x_1,...,x_{2n})=f_1(x_1,...,x_{2n})\,+\, S_{2n}  $$ (an $LPI$ with some words having sum of exponents zero at every variable).

By the same arguments we can see that if $A= M_n(R)$ be $n>1$  and $R$ a infinite commutative domain,   there is not exist a non-zero polynomial $g\, \in \, R[X] $ not zero such that $g(ab)=0$ for all $a,\, b \in \, A $ with $a^2=b^2=0 $. Note that the other side of Proposition 1.2 is true. In fact, suppose that $R$ is finite, so $A$ is finite. From this for all  matrices  $x \in A$ there exist integers $r>t>0$ such that $x^r=x^t $. So $x$ satisfies the polynomial $p_x(X)=X^r-X^t$. Then defining $g= \displaystyle\prod_{x\in A} p_x$ one has $g(ab)=0$ for all $a,\, b \in \, A $ with $a^2=b^2=0 $, which is a contradiction.  

\begin{lem1}\label{2}
	
	Let $A= M_n(R)$ be $n>1$  and $R$ a  commutative domain. Then $R$ is infinite if and only if there is not any non-zero polynomial $g\, \in \, R[X] $  such that $g(ab)=0$ for all $a,\, b \in \, A $ with $a^2=b^2=0 $.
	
\end{lem1}

\begin{prop}\label{p2}
	Let $F=\displaystyle\frac{R[x, y]}{(x^2,\, y^2)}$  where $R$ is a infinite commutative domain. Then $F$ satisfies any polynomial identity $f$ of $M_n(R)$ where $ n\geq2$. Also,  $U(F)$ does not satisfy the standard polynomial  $S_2$ and $S_3$. In particular, $F$ satisfies $S_{2n}$ if only if $n\geq 2$. 
	
	\begin{proof}
		
	Suppose by contradiction that $F$ does not satisfy an identity $f$ of $M_n(R)$ where $ n\geq2$. Then there exist   $\overline{u_1}=\overline{ g_1(x,y)},...,\overline{u_{2n}}=\overline{g_{2n}(x,y)} \, \in \, F$  with 
		$$ f(\overline{u_1},...,\overline{u_{2n}})\neq 0 .$$ So (multiplying by $x$ or $y$ on left or right if necessary) we obtain  $$g(xy) = f(u_1,...,u_{2n})\neq \sum_{i}f_iu^2g_i, $$ $u=x,\, y$. Then, for all   $a,b\, \in \, A= M_n(R)$ such that $a^2=b^2=0$, we have  $ g(ab) = 0 $  it is a contradiction by  Lemma \ref{2}. By evaluation of $S_2$ and $S_3$ on  units $\overline {1+x},\, \overline { 1 + y}$ and $\overline{(1 +x)(1+y)}$, we see that $U(F)$ does not satisfy $S_2$ and $S_3.$ 		
		
	\end{proof}
	
\end{prop}

\begin{cor1}
	
	Let $B$	be an algebra over a commutative domain $R$ whose $U(A)$ satisfies the standard polynomial  $S_3$. Then there exists a polynomial $g \in R[X]$ such that $g(ab)=0$ for all $a,\,b \in B$ with $a^2=b^2=0.$
	
\end{cor1}

\begin{proof}
	
By applying  $S_3$ on  $X, Y$ and $XY$ we obtain $$S_3(X,\, Y, \, XY)= (YX)^2 -X^2Y^2 - YX^2Y + XY^2X.$$ So $S_3(1 +x, 1+y, (1+x)(1+y))\neq 0 $ in $F$, then  multiplying by $x$ on the left and  $y$ on the right  we obtain a polynomial $g(xy)=xS_3y$ such that $g(ab)= 0.$ 	
	
\end{proof}	

\begin{cor1}
	
	Let  $F$ be the algebra in Proposition \ref{p2}. The $F$ does not satisfy the group identity $w=1$, with $R$ a infinite field. 
\end{cor1}	

\begin{proof}
	
	If $F$ satisfies a group identity $w=1$. By Lemma 3.1 \cite{LIU} there exists a polynomial $g(t) \in R[t] $ such that $g(\overline{x}\,\overline{y})=0$. Thus $\overline{g(xy)}=0$, that is, $g(xy)=\sum f_iu^2g_i, \, u=x,\, y$. Then, $g(ab)=0$ for all $a,b \in M_n(R)$ such that $a^2=b^2=0$ and this contradicts Lemma 1.1.
	
\end{proof}

With some adaptations Lemma 3.1 \cite{LIU} holds whenever $R$ is  an infinite commutative domain. It follows that Corollary \ref{c2} holds if $R$ is an infinite commutative domain. 
Note that $F$ is an algebra with $LPI$ ( having more than one word which  not constant) with $\sum \, exp w_i =0$  without group identity ($LPI$ with only one non-constant word).

\newpage

\printindex


\begin{thebibliography}{99}
	
	\bibitem{GSV} A. Giambruno, S.Sehgal, and A. Valenti, \emph{Group algebras whose units satisfy a group identity}, Proc. Amer. Math. Soc. 125 (1997), 629-634. MR 97f:16056
	
	
	\bibitem{GJV94} A. Giambruno, E. Jespers, and A. Valenti, \emph{Group identities on units of rings}, Arch. Mathj. (Basel)
	63 (1994), no. 4, 291-296. MR 95:16044
	
	\bibitem{P} D.S. Passman, \emph{The algebraic structure of group rings}, Dover, New York 2011
	
	\bibitem{GonD} J. Z. Gonçalves, M.A. Dokuchaev \emph{Identities on units of alge
		braic algebras}, Journal of Algebra 250, 638-646 (2002)
	
	\bibitem{Gon84} J.Z. Gonçalves, \emph{Free subgroups of units in group rings}, Canad. Math. Bull.27 (1984), no. 3, 309-312. MR 85k:20021
	\bibitem{LIU}Liu, Chia-Hsin, \emph{Group Algebras With Units Satisfying a Group Identity}, Proc. Amer. Math. Soc. 127, (1999f), 327-336.
	
\end{thebibliography}
\end{document}